\documentstyle[12pt,amssymb]{article}
\begin{document}
\font\frak=eufm10 scaled\magstep1
\font\fak=eufm10 scaled\magstep2
\font\fk=eufm10 scaled\magstep3
\font\scriptfrak=eufm10
\font\tenfrak=eufm10

\font\tengoth=eufm10 scaled\magstep1 \font\sevengoth=eufm7 \font\fivegoth=eufm5
\newfam\gothfam
\textfont\gothfam=\tengoth\scriptfont\gothfam=\sevengoth
  \scriptscriptfont\gothfam=\fivegoth
\def\goth{\fam\gothfam}    

\newtheorem{theorem}{Theorem}
\newtheorem{corollary}{Corollary}
\newtheorem{proposition}{Proposition}
\newtheorem{definition}{Definition}
\newtheorem{lemma}{Lemma}
\font\frak=eufm10 scaled\magstep1

\newtheorem{Theo}{Theorem}
\newtheorem{ex}{Example}
\newtheorem{cor}{Corollary}
\newenvironment{pf}{{\noindent{\it Proof.-}}}{\ $\Box$\medskip}


\mathchardef\za="710B  
\mathchardef\zb="710C  
\mathchardef\zg="710D  
\mathchardef\zd="710E  
\mathchardef\zve="710F 
\mathchardef\zz="7110  
\mathchardef\zh="7111  
\mathchardef\zvy="7112 
\mathchardef\zi="7113  
\mathchardef\zk="7114  
\mathchardef\zl="7115  
\mathchardef\zm="7116  
\mathchardef\zn="7117  
\mathchardef\zx="7118  
\mathchardef\zp="7119  
\mathchardef\zr="711A  
\mathchardef\zs="711B  
\mathchardef\zt="711C  
\mathchardef\zu="711D  
\mathchardef\zvf="711E 
\mathchardef\zq="711F  
\mathchardef\zc="7120  
\mathchardef\zw="7121  
\mathchardef\ze="7122  
\mathchardef\zy="7123  
\mathchardef\zf="7124  
\mathchardef\zvr="7125 
\mathchardef\zvs="7126 
\mathchardef\zf="7127  
\mathchardef\zG="7000  
\mathchardef\zD="7001  
\mathchardef\zY="7002  
\mathchardef\zL="7003  
\mathchardef\zX="7004  
\mathchardef\zP="7005  
\mathchardef\zS="7006  
\mathchardef\zU="7007  
\mathchardef\zF="7008  
\mathchardef\zW="700A  

\def\N{{\frak N}}
\def\const{\hbox{const}}
\def\id{\mathop{\rm id}\nolimits}
\def\<#1>{\langle#1\rangle}        

\newcommand{\be}{\begin{equation}}
\newcommand{\ee}{\end{equation}}
\newcommand{\ra}{\rightarrow}
\newcommand{\lra}{\longrightarrow}
\newcommand{\bea}{\begin{eqnarray}}
\newcommand{\eea}{\end{eqnarray}}
\newcommand{\beas}{\begin{eqnarray*}}
\newcommand{\eeas}{\end{eqnarray*}}
\newcommand{\Z}{{\Bbb Z}}
\newcommand{\R}{{\Bbb R}}
\newcommand{\C}{{\Bbb C}}
\newcommand{\1}{{\bold 1}}
\newcommand{\SL}{SL(2,\C)}
\newcommand{\Sl}{sl(2,\C)}
\newcommand{\SU}{SU(2)}
\newcommand{\su}{{\goth{su}}(2)}
\newcommand{\SB}{SB(2,\C)}
\newcommand{\Sb}{sb(2,\C)}
\newcommand{\G}{{\goth g}}
\newcommand{\D}{{\rm d}}
\newcommand{\de}{\,{\stackrel{\rm def}{=}}\,}
\newcommand{\we}{\wedge}
\newcommand{\nn}{\nonumber}
\newcommand{\ot}{\otimes}
\newcommand{\s}{{\textstyle *}}
\newcommand{\ts}{T^\s}
\newcommand{\da}{\dagger}
\newcommand{\pa}{\partial}
\newcommand{\ti}{\times}
\newcommand{\A}{{\cal A}}
\newcommand{\Li}{{\cal L}}
\newcommand{\ka}{{\Bbb K}}
\newcommand{\find}{\mid}
\def\ep{\epsilon}
\def\uno{{\,I}}

\title{Contractions: Nijenhuis and Saletan tensors for general
algebraic structures}

\author{Jos\'e F. Cari\~nena\\
Departamento de F\'{\i}sica Te\'orica, Univ. de Zaragoza\\
50.009 Zaragoza, Spain\\
{\it e-mail:} jfc@posta.unizar.es
\and
Janusz Grabowski\thanks{Supported by KBN, grant No. 2 P03A 031 17.}\\
Institute of Mathematics, Warsaw University\\
ul. Banacha 2, 02-097 Warszawa, Poland. \\
and\\
Mathematical Institute, Polish Academy of Sciences\\
ul. \'Sniadeckich 8, P. O. Box 137, 00-950 Warszawa, Poland\\
{\it e-mail:} jagrab@mimuw.edu.pl
\and
Giuseppe Marmo\thanks{Supported by PRIN SINTESI.}\\
Dipartimento di Scienze Fisiche,
Universit\`a Federico II di Napoli\\
and\\
INFN, Sezione di Napoli\\
Complesso Universitario di Monte Sant'Angelo\\
Via Cintia, 80126 Napoli, Italy\\
{\it e-mail:} marmo@na.infn.it}
\maketitle
\begin{abstract}
In this note we study  generalizations in many directions  of
the contraction procedure for Lie algebras  introduced  by  Saletan
\cite{Sa}. We consider products  of  arbitrary
nature,  not  necessarily  Lie brackets, and we generalize  to
infinite  dimension, considering  a modification  of  the  approach  by
Nijenhuis  tensors   to   bilinear operations on sections of
finite-dimensional  vector  bundles. We apply our general procedure to Lie
algebras, Lie algebroids, and Poisson brackets. We present also results on
contractions of $n$-ary products and coproducts.
\end{abstract}

\section{Introduction}

For a general (real) topological algebra,  i.e.,  a  topological  vector
space $\A$ over $\R$ (but other topological fields,  like  $\C$,  can  be
considered in  a  similar  way  as  well)  with  a  continuous  bilinear
operation
\be
\mu:\A\times\A\rightarrow\A,\quad (X,Y)\mapsto X*Y,
\ee
one considers contraction procedures as follows.
\par
If $U(\zl):\A\ra\A$ is a family of linear morphisms  which  continuously
depends on the parameter $\zl\in\R$ from a neighbourhood $\cal U$  of  0
and $U(\zl)$ are invertible for $\zl\in{\cal U}\setminus\{ 0\}$, then we
can consider the continuous family of products $X*^\zl Y$ defined by
\be
X*^\zl Y=U(\zl)^{-1}(U(\zl)(X)*U(\zl)(Y)),
\ee
for $\zl\in{\cal U}\setminus\{ 0\}$. All these products are   isomorphic
by definition, since
\be
U(\zl)(X*^\zl Y)=U(\zl)(X)*U(\zl)(Y)
\ee
and if $N=U(0)$ is invertible, then clearly
\be\label{N}
\lim_{\zl\to 0}X*^\zl Y=N^{-1}(N(X)*N(Y)).
\ee
But sometimes, the limit $\lim_{\zl\to 0}X*^\zl Y$  may  exist  for  all
$X,Y\in\A$ even if $N$ is not invertible and  (\ref{N})  does  not  make
sense. We say then that $\lim_{\zl\to 0}X*^\zl Y$ is  a {\it contraction}
of the product $X*Y$. Of course, the problem of existence and the  form
of the contracted product  heavily  depends  on  the  family  $U(\zl)$.
In \cite{Sa} this problem has been solved for linear  families
$U(\zl)=\zl\,I+N$ and $\A$ -- a finite-dimensional Lie algebra.
\par
Here we study  generalizations in various directions  of  the  contraction
procedure introduced by Saletan \cite{Sa}. First  of
all, we consider products  of  arbitrary  nature,  not  necessarily  Lie
brackets. Second, we generalize  to  infinite  dimension, considering  a
modification  of  the  approach  by  Nijenhuis  tensors   to   bilinear
operations on sections of  finite-dimensional vector bundles.
The motivation stems from physics, since infinite-dimensional
algebras of sections of some bundles arise frequently as models both
in Classical and Quantum Physics.
In particular, we were confronted with this problem within the framework
of Quantum Bihamiltonian Systems \cite{CGM}. According to Dirac
\cite{Di}, a "quantum Poisson bracket" necessarily arises from the
associative product on the space of operators. Similarly, by Ado's
theorem, any finite-dimensional Lie algebra arises as an algebra of
matrices.
It is therefore quite natural to investigate contractions of associative
algebras along with contractions of Lie algebras and their generalizations
to Lie algebroids.
We concentrate mainly on smooth sections, but this particular choice plays
no definite role in our  approach.

The paper is organized as follows. In the next section we present the
general scheme we are working with and the main result (Theorem 2) on
contractions for algebras of sections of vector bundles. We make some
remarks on contractions with respect to  more general linear families
$U(\zl)=\zl\, A+N$.

Section 3 is devoted to examples and Section 4 to more detailed studies of
hierarchies of contractions. We develop an algebraic technique which
allows us to produce much simpler proofs of facts about hierarchies than
those available in the literature.

In Section 5 we comment on the behaviour of algebraic properties under
contractions.

Contractions of Lie algebras and Lie algebroids, as particular cases of
our general procedure, are studied in Sections 6 and 7.

In Section 8 we use our knowledge on contractions of Lie algebroids to
define contractions of Poisson structures. The approach is very natural
and leads to structures very similar (but slightly different) to those
which are known under the name of Poisson-Nijenhuis structures (cf.
\cite{MM,KSM}).

We end up with  observations on contractions of $n$-ary products and
coproducts.

\section{Linear contractions of products on sections of vector bundles}

Let us assume that $E$ is a smooth vector bundle with fibers of dimension
$n_0$,
over a smooth manifold $M$. Denote by $*$ a bilinear operation
$\mu:{\cal A}\times {\cal A}\to {\cal A}$,
$$\mu: (X,Y)\in {\cal A}\otimes {\cal A}\mapsto X*Y\in {\cal A}\ ,
$$
on the space of smooth sections of  $E$  which is, at  least point-wise,
continuous. In practice, we shall deal with local products, therefore
being defined by bilinear differential operators. We shall use both
notations for the product according to which one is more convenient when
treating  particular cases.
\par
Let $N:E\to E$ be a
smooth vector bundle morphism over $\id_M$. One refers also to $N$ as to a
(1,1)--tensor field, i.e., a section of $E^*\otimes E$. Since $N_p:E_p\to
E_p$ is a morphism of the finite-dimensional
vector space $E_p$, where $E_p$ denotes the fiber over the point $p\in M$,
 we have the Riesz decomposition $E_p=E_p^1\oplus E_p^2$ into invariant
subspaces of $N_p$ in such a way that $N_p$ is invertible on $E_p^1$ and
nilpotent of order $q$
 on $E_p^2$, i.e., $N^q(X_p)=0$, for $X_p\in E_p^2$. One can take
$E_p^1=\widetilde N_p(E_p)$, $E_p^2=\ker \widetilde
N_p$,
 where $\widetilde N_p=(N_p)^{n_0}$, with $n_0=\dim E_p$. In this way
we get the decomposition
 $E=E^1\oplus E^2$ of the vector bundle $E$ into two supplementary
generalized distributions. Note that the dimension of $E_p^1$ may vary
from point to point. Nevertheless, $E^1$ is a smooth distribution,
i.e., it is  generated locally by a finite number of smooth sections of $E$.
 Indeed, if $\{e_1,\ldots,e_{n_0}\}$ is a local basis of smooth sections
of $E$, then $\{\widetilde N(e_1),\ldots, \widetilde N(e_{n_0})\}$ is a
set of local smooth sections generating locally $E^1$.
\begin{theorem}
The (generalized) distribution $E^2$ is smooth if  and  only  if  it  is
regular, i.e., of constant rank: $\dim E_p^2=\const.$
\end{theorem}
\begin{pf} Since the rank of a smooth distribution is semi-continuous from
above:
\begin{equation}
\lim_{p\to p_0}\inf\dim E_p^2\geq \dim E_{p_0}^2\ ,\label{firstc}
\end{equation}
and the complementary distribution $E^1$ is smooth, so that
\begin{equation}
\lim_{p\to p_0}\sup\dim E_p^2\leq \dim E_{p_0}^2\ ,\label{secondc}
\end{equation}
we see that both conditions (\ref{firstc}) and (\ref{secondc})
are satisfied if and only if $E^2$ is of constant rank.

Conversely, if  $E^2$ is of constant rank, say $n_0-l$, take a basis
 $\{e_1,\ldots,e_{n_0}\}$ of smooth local sections of $E$ such that
the elements of  $\{\widetilde N(e_1),\ldots, \widetilde N(e_{l})\}$
span $E_p^1$. Then $\{\widetilde N(e_1),\ldots, \widetilde N(e_{l})\}$
is a basis of local sections of $E^1$ near $p\in M$. Write
$$\widetilde N(e_i)=\sum_{j=1}^lf_{ij}\, \widetilde N(e_j)\ .
$$  Then the functions $f_{ij}$ are smooth and the smooth sections
$$\widetilde e_i=e_i-\sum_{j=1}^l f_{ij}\, e_j\ ,\qquad i=l+1,\ldots,n_0\,,
$$
span locally $E^2$. Indeed, $\widetilde N(\widetilde e_i)=0$ and
the elements $\widetilde e_i$, for $i=l+1,\ldots,n_0$,
are linearly independent.
\end{pf}

Note that in general no one of the distributions $E^1$ and $E^2$ has
to be \lq\lq involutive'' in the sense that smooth sections of $E^1$
(resp. $E^2$) are closed with respect to the composition law $*$.

Consider now a new (1,1)-tensor $U(\lambda)=\lambda \, I+N$ depending on a
real parameter $\lambda$. Since the spectrum of $N$ is finite and
continuously depends on $p$, in a sufficiently small neighbourhood of $p$
all of $U(\lambda)_p$ are invertible for sufficiently small $\lambda$, but
$\lambda\ne 0$.
Thus, we can locally define, for $\lambda\ne 0$, a new operation
\begin{eqnarray}
X*_N^\zl Y&=&U(\lambda)^{-1}(U(\lambda)(X)*U(\lambda)(Y))\cr
&=&U(\lambda)^{-1}((\lambda\, X+N(X))*(\lambda\, Y+N(Y)))\cr
&=&U(\lambda)^{-1}(\lambda^2\, X*Y+\lambda\,(N(X)*Y+X*N(Y))+N(X)*N(Y))\
.\label{starlam}
\end{eqnarray}
We would like to find conditions assuring that the limit
$$X*_NY=\lim_{\lambda\to 0} X*_N^\lambda Y$$
exists for all $X,Y\in {\cal A}$ and find the corresponding contraction
$X*_NY$.

Using the identity $U(\lambda)^{-1}(\lambda\, I+N)=I$, i.e.,
$$U(\lambda)^{-1}(\lambda \, X)=X-U(\lambda)^{-1}N(X)\ ,
$$
we get from (\ref{starlam}) that
\begin{eqnarray}
X*_N^\lambda Y&=&\lambda\, X*Y+(N(X)* Y+X*N(Y)-N(X*Y))\cr
&+&U(\lambda)^{-1}(N(X)*N(Y)-N(N(X)*Y+X*N(Y)-N(X*Y)))\ .\label{starlam2}
\end{eqnarray}
Denoting
$$\delta_N\mu(X,Y)=X\widetilde *_NY=N(X)*Y+X*N(Y)-N(X*Y)\ ,
$$
and by $T_N\zm(X,Y)$ -- the {\it Nijenhuis torsion} of N:
$$T_N\zm(X,Y)=N(X)*N(Y)-N(X\widetilde *_NY)\ ,
$$
we can rewrite (\ref{starlam2}) in the form
$$X*_N^\lambda Y=\lambda\, X*Y+X\widetilde*_NY+U(\lambda)^{-1}T_N\mu(X,Y)\ .
$$
Hence, the limit  $$\lim_{\lambda\to 0} X*_N^\lambda Y$$ exists  if  and
only if
\be
\lim_{\lambda\to 0} U(\lambda)^{-1} T_N\mu(X,Y)  \qquad {\rm exists\ for\
every}\quad X,Y\in
{\cal A}\ . \label{limcond}
\ee

Denote by ${\cal A}^1$, ${\cal A}^2$, the spaces of smooth sections of
$E^1$ and $E^2$, respectively. Of course, in general ${\cal A}\ne{\cal
A}^1\oplus {\cal A}^2$. We may  have  ${\cal A}^2=\{  0\}$  even  in  the
case $E^2\ne\{ 0\}$.
Since $E^1$ and $E^2$ are invariant distributions of $U(\lambda)$, hence
of $U(\lambda)^{-1}$, the existence of the  limit  (\ref{limcond})  may
be checked separately on the corresponding parts of $T_N\mu$. On $E^2$ the
tensor $N$ is nilpotent, so  for $X_p\in E_p^2$,
\begin{equation}
(\lambda\, I+N)_p^{-1}(X_p)=\left(\lambda(I-(-N/\zl))\right)_p^{-1}(X_p)=
\frac 1{\lambda}\sum_{n=0}^\infty\left(-\frac 1\lambda\right)^nN_p^n(X_p)\ ,\label{expinv}
\end{equation}
where the sum is in fact finite, and
$$\lim_{\lambda\to
0}\sum_{n=0}^{q-1}\frac{(-1)^n}{\lambda^{n+1}}\,N_p^n(X_p)
$$
exists if and only if $X_p=0$. Thus, a necessary condition for existence
of the limit (\ref{limcond}) is that $T_N\mu(X,Y)\in {\cal A}^1$ for every
$X,Y\in {\cal A}$.

Since on $E^1$ the tensor $N$ is invertible, we have clearly
$$\lim_{\lambda\to 0}(\lambda\, I+N)^{-1}=N^{-1}$$
on $E^1$, so that, assuming  $T_N\mu(X,Y)\in {\cal A}^1$,
$$\lim_{\lambda\to 0}U(\lambda)^{-1} T_N\mu(X,Y) =N^{-1} T_N\mu(X,Y)
=\tau_N\mu(X,Y)\ .
$$
Here $\tau_N\mu(X,Y)=N^{-1} T_N\mu(X,Y)$ is the unique section of $E^1$
determined by the condition
$$N(\tau_N\mu(X,Y))=T_N\mu(X,Y)\ .
$$
In order to get a new product on ${\cal A}$ we have to assume that
$\tau_N\mu(X,Y)$ is smooth, which is {\it a priori} not automatic,
even if we have $T_N\mu(X,Y)\in  {\cal A}^1$.
Note that if $N$ is regular, i.e., $E^1$ is of constant dimension, then,
as we shall show in Theorem 3, $N({\cal A}^1)= {\cal A}^1$ and
$\tau_N\mu(X,Y)$ is smooth automatically.

Let us summarize the above as follows:
\begin{theorem}
Let  $\mu:{\cal  A}\times  {\cal  A}\to  {\cal  A}$   be   a   point-wise
continuous bilinear product of smooth sections of a vector bundle $E$
over a manifold $M$
(we will write also $X*Y$ instead of $\mu(X,Y)$) and let $N:E\to E$ be a
smooth (1,1)--tensor.
Denote by $U(\lambda)=\lambda \, I+N$ a deformation of $N$, by
$E=E^1\oplus E^2$ the Riesz decomposition of $E$ relative to $N$, and by
${\cal A}^1$ the set of smooth sections of $E^1$.
Then, the limit
$$\lim_{\lambda\to 0}U(\lambda)^{-1}(U(\zl)(X)*U(\zl)(Y))
$$
exists for all $X,Y\in\A$ and defines a new (contracted) bilinear
operation
$$D_N\zm(X,Y)=X*_NY$$
on $\cal A$ if and only if the Nijenhuis torsion
$$T_N\mu(X,Y)=N(X)*N(Y)-N(N(X)*Y+X*N(Y))+N^2(X*Y)
$$
takes values in $N({\cal A}^1)$. If this is the case, then
\be
X*_NY=X\widetilde *_NY+\tau_N\mu(X,Y)\ ,\label{tildest}
\ee
where $X\widetilde *_NY$ is a new  bilinear  operation $\delta_N\mu$  on
$\A$ defined by
$$\zd_N\zm=X\widetilde *_NY=N(X)*Y+X*N(Y)-N(X*Y)\ ,
$$
and  $\tau_N\mu(X,Y)=N^{-1} T_N\mu(X,Y)$ is the unique section of ${\cal
A}^1$
such that $$N(\tau_N\mu(X,Y))=T_N\mu(X,Y)\ .
$$
Moreover, $N$ constitutes a homomorphism of $({\cal A},\mu_N)$ into
$({\cal A},\mu)$:
$$N(X*_NY)=N(X)*N(Y)\ .
$$
\end{theorem}
{\bf Remark.}  Let us note that our procedure is not just applying the
finite--dimensional linear one to every fiber, since the operation $*$
need not act fiber-wise. Also, this is not direct application to
infinite--dimensional algebra ${\cal A}$, since we have not,
in general, the  Riesz  decomposition  ${\cal  A}={\cal  A}_1\oplus{\cal
A}_2$ with respect to $N$. On the other hand, the whole procedure can be
applied directly to infinite-dimensional cases for which we are given
the Riesz decomposition of $N$.

\begin{definition} The tensor $N$ satisfying the assumptions of Theorem 2,
i.e., such that $T_N\zm$ takes values in $N(\A^1)$, will be called {\sl a
Saletan tensor}. If $N^k$ is a Saletan tensor for every $k=1,2,3,\dots$,
then $N$ will be called a {\sl strong Saletan tensor}. In the case
$T_N\zm=0$ we shall call $N$ a {\sl Nijenhuis tensor}.
\end{definition}
{\bf Remark.} It is obvious from the proof of Theorem 2 that Nijenhuis
tensors define contractions even in the case of infinite-dimensional
algebras $\A$ without any assumption that $\A$ consists of sections of a
finite-dimensional vector bundle. Indeed, with $T_N\zm=0$ we have no
obstructions, the Riesz decomposition is irrelevant, and
$X*_NY=X\tilde*_NY$.
Of course,

(Nijenhuis) $\Rightarrow$ (strong Saletan) $\Rightarrow$ (Saletan).

\noindent
We shall call $N$ {\it regular}, if  $E^1$
(hence also $E^2$) is of constant rank. This is always the case when
$E$ is a bundle over a single point, i.e., $E=\A$.

\begin{theorem} In the regular case, i.e., when $E^1$ is of constant rank,
$N({\cal A}^1)={\cal A}^1$, so that $N$ is a Saletan tensor if and only if
$T_N\mu(X,Y)$  takes values in $\A^1$.
\end{theorem}
{\bf Remark.} We shall prove a stronger result in Theorem 6.

\begin{pf}
Indeed, according to Theorem 1, both $E^1$ and $E^2$ are smooth
distributions. Locally we have a basis of smooth sections of $E^1$, and
$N$ acts on this basis as invertible matrix of smooth functions.
Indeed, since regularity of $E^1$ implies that there is a local basis
$\{e_1,\ldots,e_l\}$ of sections of $E^1$, in this basis $N$ acts simply
as invertible matrix of smooth functions $(f_{ij})$, so for $X=\sum g_i\, e_i$,
$$N \left(\sum_{i=1}^l h_i\, e_i\right)=X,
$$
where the smooth functions
$h_i$ are defined by
$$\sum_{j=1}^l f_{ij}\, h_j=g_i\ .
$$
Hence, $N^{-1}(X)$ is locally, thus globally, smooth section of $E^1$ for
any smooth section $X$ of $E^1$.
\end{pf}

\noindent
{\bf Remark.} For a fiber bundle over a single point Theorem
2 gives exactly the Saletan result \cite{Sa} in case $\mu$ is a Lie
bracket. Saletan writes $X*_NY$ in the form
\be\label{saldef}
X*_NY=(X\widetilde *_NY)_2+N^{-1}\left((N(X)*N(Y))_1\right)\ ,
\ee
where $X=X_1+X_2$ is the decomposition of $X\in {\cal A}$
into sections of $E^1$ and  $E^2$.
Of course, (\ref{saldef}) is formally the same as (\ref{tildest}) for the
decomposition into sections of $E^1$ and $E^2$. However, in general the
summands of the right hand side of (\ref{saldef}) are not smooth, while
the decomposition (\ref{tildest}) is into smooth parts. In  the  regular
case both  formulae coincide.

\begin{theorem}
(a)  Theorem  2  remains  valid  when we
consider  the  family  $U(\zl)$  in  a  slightly  more   general   form:
$U(\zl)=\zl \,I+f(\zl)\,N$, where $f$ is  continuous  and  $f(0)=1$.

(b) If   we   consider    instead of $U(\zl)$    the    family
$U_1(\zl)=\zl A+N$, then
the contraction procedure for $U_1(\zl)$ and the product  $*$  is
equivalent to the contraction procedure of the above type for a new $N$
and a new product. In particular, if $A$ is invertible, we get our
standard  contraction for $A^{-1}N$  and the product
$X*_AY=A^{-1}(A(X)*A(Y))$. In other words, the contraction procedure for
the family $U_1(\zl)=\zl\, A+N$ can be reduced to the contraction
described in Theorem 2.
\end{theorem}
\begin{pf}

(a) Let us write $U(\zl)=\frac{U_1(\ze)}{f(\zl)}$, where
$U_1(\ze)=\ze\, I+N$ and $\ze=\frac{\zl}{f(\zl)}$, so that $\zl\ra 0$ is
equivalent to $\ze\ra 0$. Since
\be
U(\zl)^{-1}(U(\zl)(X)*U(\zl)(Y))=\frac{1}{f(\zl)}\,
U_1(\ze)^{-1}(U_1(\ze)(X)*U_1(\ze)(Y))
\ee
and $\lim_{\zl\to 0}f(\zl)=1$,
both contraction procedures are equivalent.

(b) Assume first that $A$ is invertible.
Since $\zl\, A+N=A(\zl\, I+A^{-1}N)$, we can use Theorem 2 for
$N:=A^{-1}N$ and the product $*_A$.  In  fact,  we  can  skip  the
assumption  that  $A$ is invertible. For, take $\zl_0$  for  which
$A+\zl_0\N$  is  invertible. Then, we write
\be
U(\zl)=\zl\, A+N=(A+\zl_0\, N)\, (\zl\,  I+(1-\zl\,\zl_0)(A+\zl_0\, N)^{-1}N)\ ,
\ee
and we can proceed as before and using (a) of the Theorem.
\end{pf}

\section{Examples}

Many interesting physical applications are based on the idea of
contraction by \.In\"on\"u and Wigner \cite{IW}. We will call a smooth
distribution $E^1$ in the vector bundle $E$ {\it involutive} if the space
$\A^1$ of sections of $E^1$ is closed with respect to the product $*$,
i.e., $\A^1$ is a subalgebra of $\A$.
\begin{theorem} Let $E^1$ be a smooth regular and
involutive distribution in $E$. Take $E^2$ to be any supplementary
smooth distribution and let $N=P_{E^1}$ be the projection on $E^1$ along
$E^2$.
Then $N$ is a Saletan tensor which is Nijenhuis if and only if $E^2$ is
also involutive. The contracted product reads
\be
X*_NY=X_1*Y_1+(X_1*Y_2+X_2*Y_1)_2,
\ee
where $X=X_1+X_2$, etc.,
is the decomposition with respect to the splitting $E=E^1\oplus E^2$.
\end{theorem}
\begin{pf}
It is obvious that the Nijenhuis tensor $T_N(X,Y)=N(X)*N(Y)-N(X\widetilde
*_NY)$ takes values in $\A^1$, since $E^1$ is involutive. Due to
regularity, the corresponding contraction exists (Theorem 3). It is easy
to see that
\be
X\widetilde*_NY=X_1*Y_1+(X_1*Y_2+X_2*Y_1)_2-(X_2*Y_2)_1\, .
\ee
Hence, $T_N(X,Y)=\tau_N(X,Y)=(X_2*Y_2)_1$, so that $N$ is a Nijenhuis
tensor if and only if $E^2$ is also involutive.  Finally,
\be
X*_NY=X\widetilde *_NY+\tau_N(X,Y)=X_1*Y_1+(X_1*Y_2+X_2*Y_1)_2\, .
\ee
\end{pf}

\medskip\noindent
{\bf Example 1.} Consider a manifold $M$ with two foliations ${\cal F}_1,
{\cal F}_2$ corresponding to a splitting into complementary distributions
$TM=E^1\oplus E^2$. The projection $N$ of $TM$ onto $E^1$ along $E^2$ is a
Nijenhuis tensor (Theorem 5). The contracted bracket is trivial for
two vector fields which are tangent to ${\cal F}_2$, it is the standard
one for two vector fields which are tangent to ${\cal F}_1$ and it is the
projection onto $E^2$ of the standard bracket of two vector fields of
which one belongs to ${\cal F}_1$ and the second to ${\cal F}_2$.

\medskip\noindent
{\bf Example 2.} Let $E$ be  just  1--dimensional  trivial  bundle  over
$\R$, i.e., ${\cal A}=C^\infty(\R)$.
Take $f*g=f'\, g'$ and $N=\varphi \, I$, where $\varphi \in{\cal A}$. Then
$E^1_p=T_p\R$
if $\varphi(p)\ne 0$ and $E^1_p=\{0\}$ otherwise, so the distribution
need not to be regular.
We have
$$f\widetilde *_Ng=\varphi\, f'\,g'+\varphi'\, (f'\,g+f\,g')$$
and $T_N\mu(f,g)=\varphi'\, \varphi'\, f\, g$.
For instance, if $\varphi(p)=p^2$ (non--regular case), then
$$
T_N\mu(f,g)=4\varphi \,f\,g\, ,$$
i.e.,
$N$ is not Nijenhuis but satisfies the assumptions of the Theorem. We get
$$
f*_Ng=f\widetilde *_Ng+N^{-1}(T_N\mu(f,g))=\varphi\, f'\,g'+\varphi'
(f'g+fg')+4f\, g\ .
$$

\medskip\noindent
{\bf Example 3.} It is easy to see that if $*$ is an associative product,
the multiplication by any $K\in\A$:
$$N_K:\A\ra\A,\quad N_K(X)=KX,
$$
is  a  Nijenhuis  tensor.  In  view  of  Remark  4,  the   corresponding
contraction yields
$$X*_{N_K}Y=X*K*Y.
$$
This product has been recently used as an alternative product of
operators in Quantum Mechanics in connection with deformed oscillators
\cite{MMSZ}, taking up an old idea of Wigner \cite{Wi}.

\medskip\noindent
{\bf Example 4.} Another alternative product for Quantum Mechanics can be
constructed as a contraction as follows (cf. \cite{CGM}).
Let  now  the  algebra  $\A$  be  the  algebra  of $n\ti n$ matrices,
$n=1,2,\dots,\infty$. In the case $n=\infty$ we consider infinite matrices
concentrated on the diagonal, i.e., matrices which are  null  outside  a
strip of  the  diagonal.  The  algebra  $\A$  represents  then  unbounded
operators on a Hilbert space $\cal H$ with a  common  dense  domain.  We
choose $\A_1$ to be a subalgebra of upper-triangular  matrices  and  for
$\A_2$ we take the supplementary  algebra  of  strict  lower-triangular
matrices.  Then,  the   mapping
\be
N_\za(A)=(1-\za)A_1+\za A
\ee
is a Nijenhuis tensor  on  $\A$ for  every  $\za\in{\Bbb  C}$.
For example, for $n=2$, the  new associative matrix multiplication has
the form
\be
\left(\begin{array}{cc}a&b\\c&d\end{array}\right)\circ
\left(\begin{array}{cc}a'&b'\\c'&d'\end{array}\right)=
\left(\begin{array}{cc}aa'+\za
bc'&ab'+bd'\\ca'+dc'&dd'+\za cb'\end{array}\right).
\ee
Note that the unit matrix $I$ remains the unit for this new product  and
that inner derivations given by diagonal matrices are the same for  both
products.

Since  the  corresponding   deformed
associative products $*_{N_\za}$ give all the same result if one  of
factors is a diagonal matrix, in the infinite case $n=\infty$  the
Hamiltonian $H$ for the harmonic
oscillator, $H\mid e_n\rangle=n\mid e_n\rangle$, describes the same
motion for all deformed brackets. This time, however, $a^\da*_{N_\za}
a=\za H$, so $a^\da$ and $a$ commute for $\za=0$.

\section{Hierarchy of contractions}

Let us have a look at the process of constructing contracted products in a
more systematic
way. For, denote the linear space of all bilinear products on ${\cal A}$ by
${\cal B}$,
the linear subspace of all bilinear products $\mu$ such that $\mu(X,Y)\in
N^k({\cal A}^1)$
by ${\cal B}_k^1$. Note that the distribution $E^1$ associated with the
(1,1)--tensor field $N$(N will be fixed) is the same
for all positive powers of $N$. Let
$A_N, B_N, C_N:{\cal B}\to {\cal B}$ be given by
\bea
(A_N\mu)(X,Y)&=&N(\mu(X,Y))\ ,\\
(B_N\mu)(X,Y)&=&\mu(N(X),Y)\ ,\\
(C_N\mu)(X,Y)&=&\mu(X,N(Y))\ .
\eea

It is easy to see that $A_N, B_N, C_N$ generate a commutative algebra of
linear operators on
${\cal B}$ for which ${\cal B}_k^1$ are invariant subspaces. Moreover,
$A_{N^k}=(A_N)^k$, etc.
Observe that for the derived product,
$$(\delta_N\mu)(X,Y)=\mu(N(X),Y)+\mu(X,N(Y))-N(\mu(X,Y))\ ,
$$
we can write
$$\delta_N=B_N+C_N-A_N\ ,
$$
and for the Nijenhuis torsion,
$$
(T_N\mu)(X,Y)=\mu(N(X),N(Y))-N(\delta_N\mu(X,Y))\ ,
$$
we can write
$$T_N=B_NC_N-A_N\delta_N=(A_N-B_N)(A_N-C_N)\ .
$$
The contracted product $D_N\mu$ is defined via the formula
$$\zm_N=D_N\mu= \delta_N\mu+\tau_N\mu\ ,
$$
where  $\tau_N\mu\in {\cal B}^1$ is such that $A_N\tau_N\mu=T_N\mu\ .$
Hence
\be
A_ND_N\zm=(A_N\zd_N+T_N)\zm=B_NC_N\zm.
\ee
If we use $N^k$ instead of $N$, we can define the corresponding contracted
product $D_{N^k}\mu$ if only $T_{N^k}\mu\in {\cal B}^1_k$. If this is the
case, we call such (1,1)--tensor field $N$ a {\it strong Saletan tensor}
(for $\mu$). We have the following:
\begin{theorem} If  $N$ is regular (e.g. $E$ is over a single point)
and $T_N\zm$ takes values in $\A^1$ (i.e., $N$ is a Saletan tensor),  then
$N$ is a strong Saletan tensor.
\end{theorem}
\begin{pf}
Indeed, in this case,
$$T_{N^k}\mu=(A^k_N-B^k_N)(A^k_N-C^k_N) =\omega(A_N,B_N,C_N)\,
(A_N-B_N)(A_N-C_N)\mu\ ,
$$
where $\omega$ is a polynomial and $(A_N-B_N)(A_N-C_N)\mu=T_N\zm\in{\cal
B}^1_0$, since $N$ is a Saletan tensor. We have then $T_{N^k}\mu\in{\cal
B}^1_0$, since ${\cal B}^1_0$ is an invariant subspace with respect to
$A_N, B_N, C_N$. But in the regular case $N^k(\A^1)=\A^1$ (Theorem 3),
so ${\cal B}^1_0={\cal B}^1_k$.
\end{pf}

There is a nice algebraic condition which assures that the tensor is
regular.
\begin{theorem} Suppose that there is a finite-dimensional $N$-invariant
subspace $V$ in $\A$ which generates $\A$ as a $C^\infty(M)-module$, i.e.,
the sections from $V$ span the bundle $E$.
Then the tensor $N$ is regular and it is a strong Saletan tensor if and
only if its Nijenhuis torsion takes values in $\A^1$.
\end{theorem}
\begin{pf} Let $V=V_1\oplus V_2$ be the Riesz decomposition of $V$ with
respect to $N$ (as acting on $V$). Since $N^k(V_2)=\{ 0\}$ for a
sufficiently large $k$, $V_2\subset\A_2$. Similarly, since $N(V_1)=V_1$,
$V_1\subset\A_1$. Since $V$ generates $E$, we have the
decomposition $E(p)=V_1(p)\oplus V_2(p)$ for
any $p\in M$. By the dimension argument, $V_2(p)=E_2(p)$ for every $p\in
M$, so $E^2$ is a smooth distribution and its dimension is constant due to
Theorem 1.
\end{pf}

\noindent
For Nijenhuis tensors we have the following.
\begin{theorem} If $N$ is a Nijenhuis tensor for the product  $\zm$  and
$w,v$  are  polynomials,  then  $w(N)$  is  a   Nijenhuis   tensor   for
$\zd_{v(N)}\zm$.
\end{theorem}
\begin{pf}
N a is Nijenhuis tensor  for  $\zm$, so   $T_N\zm=0$.   Since
$A_{w(N)}=w(A_N)$, etc., we have
$$T_{w(N)}\zd_{v(N)}\zm=\zd_{v(N)}W(A_N,B_N,C_N)\zm,
$$
where
$$W(x,y,z)=(w(x)-w(y))(w(x)-w(z))=W_1(x,y,z)(x-y)(x-z)
$$
for certain polynomial $W_1$. Hence,
$$T_{w(N)}\zd_{v(N)}\zm=\zd_{v(N)}W_1(A_N,B_N,C_N)T_N\zm=0.
$$
\end{pf}

\noindent
For any strong Saletan tensor $N$ we get a whole hierarchy of contracted
products
$$D_{N^k}\mu=\delta_{N^k}\mu+\tau_{N^k}\mu\ ,\qquad k=1,2,\dots
$$
We will show that this is exactly the same hierarchy if we apply
the contraction procedure
 inductively:
$$\mu_0=\mu\,,\qquad \mu_{k+1}=D_N\mu_k\ .
$$
For the case of Nijenhuis tensors, it is very easy. Indeed, as above,
$N^k$ are Nijenhuis tensors for $\mu$ for any $k=1,2,\ldots$ and
$D_{N^k}\mu=\delta_{N^k}\mu$.
To see that  $\delta_{N^k}\mu=(\delta_{N})^k\mu$, it  is  sufficient  to
check that
$$\left(\delta_{N^k}-(\delta_{N})^k\right)\mu=\left((B_N^k+C_N^k-A_N^k)-
(B_N+C_N-A_N)^k\right)\mu=0\ .
$$
But the polynomial $(x^k+y^k-z^k)-(x+y-z)^k$ vanishes for $x=z$ and for
$y=z$,
so that it can be written in the form $\omega(x,y,z)(z-x)(z-y)$. Hence,
$$\left(\delta_{N^k}-(\delta_{N})^k\right)\mu=\omega(A_N,B_N,C_N)(A_N-B_
N)(A_N-C_N)\mu=0\ ,
$$
since
$(A_N-B_N)(A_N-C_N)\mu=0$. For an arbitrary strong Saletan tensor the
situation is a little bit more complicated. First, we show the following:
\begin{lemma} With the previous notation, for any strong Saletan
tensor and any couple of natural numbers $i,k\in {\Bbb N}$, we have
\be
T_{N^i}D_{N^k}\mu=A^i_N(D_{N^{k+i}}\mu-\delta_{N^i}\,D_{N^k}\mu)\
.\label{funrel}
\ee
\end{lemma}
\begin{pf} First of all, let us observe that both sides belong to ${\cal
B}_0^1$.

\noindent
Indeed,
the left hand side equals
$$T_{N^i}(\delta_{N^k}\mu+\tau_{N^k}\mu)=\delta_{N^k}T_{N^i}\mu+T_{N^i}\tau_{N
^k}\mu\ ,
$$
and $T_{N^i}\mu,\, \tau_{N^k}\mu\in {\cal B}^1_0$, so the left hand side also
belongs
to ${\cal B}^1_0$, due to the invariance of ${\cal B}^1_0$.
As for the right hand side, we write
$$
A_N^i(D_{N^{k+i}}\mu-\delta_{N^i}D_{N^{k}}\mu)=A^i_N(\delta_{N^{k+i}}\mu
-\delta_{N^{i}}\delta_{N^{k}}\mu+\tau_{N^{k+i}}\mu-\delta_{N^{i}}
\tau_{N^{k}}\mu)\ .
$$
Since, similarly as above, $\tau_{N^{k+i}}\mu$ and $\delta_{N^{i}}
\tau_{N^{k}}\mu$ belong to ${\cal B}_0^1$, it suffices to check that
$$
\left(\delta_{N^{k+i}}-\delta_{N^{i}}\delta_{N^{k}}\right)\mu \in {\cal
B}_0^1\ ,
$$
which is straightforward, since
\bea
\left(\delta_{N^{k+i}}-\delta_{N^{i}}\delta_{N^{k}}\right)\mu=(B_N^{k+i}+C_N^{
k+i}
-A_N^{k+i}
-(B_N^i+C_N^i-A_N^i)(B_N^k+C_N^k-A_N^k))\mu\cr
=\omega(A_N,B_N,C_N)(A_N-B_N)(A_N-C_N)\mu=
\omega(A_N,B_N,C_N)T_N\mu\in {\cal B}_0^1\ ,\nonumber
\eea
where we use an analogous polynomial factor argument as above and the
invariance
of ${\cal B}_0^1$.
Hence, we can check the following by applying $A_N^k$ to  both  sides
of (\ref{funrel}) ($A_N$ is invertible on $E^1$):
 $$
 A_N^k T_{N^{i}}D_{N^{k}}\mu= A_N^{i+k}(D_{N^{i+k}}\mu
-\delta_{N^{i}}D_{N^{k}}\mu)\ .
 $$
Writing down expressions for $D_{N^{k}}\mu$ and $D_{N^{k+i}}\mu$ explicitly,
and using
$$
A_N^kD_{N^k}\mu =(A_N^k\delta_N^k+T_{N^k})\zm=B_{N}^kC_N^k\mu,
$$
etc., we get
\bea
A_N^kT_{N^i}D_{N^k}\mu&=&T_{N^i}B_N^kC_N^k\mu=(B_N^iC_N^i-\delta_{N^i}A_N^i)B_
{N}^kC_N^k\mu\\
&=&(B_N^{i+k}C_N^{i+k}-A_N^i\delta_{N^i}B_N^kC_N^k)\mu\\
&=&A_N^{i+k}(D_{N^{i+k}}\mu-\delta_{N^i}D_{N^k}\mu)\ .\nonumber
\eea
\end{pf}
\begin{cor} The tensor $N$ is a strong Saletan tensor for any of
$D_{N^k}\mu$, $k=0,1,2\ldots$.
\end{cor}

\begin{theorem} If $N$ is a strong Saletan tensor for $\mu$, then

i) We have a well--defined hierarchy of contracted products
$D_{N^k}\mu$,  $k=0,1,2\ldots$.

ii)  $N$ is a strong Saletan tensor for every $D_{N^k}\mu$,
$k=0,1,2\ldots$.

iii) $D_{N^i}D_{N^k}\mu=D_{N^{i+k}}\mu$, for any couple of natural numbers
$i,k\in {\Bbb N}$.

iv) $N^k$ is a homomorphism of the product $D_{N^{i+k}}\mu$ into $D_{N^i}\mu
$,
for any couple of natural numbers $i,k\in {\Bbb N}$.
\end{theorem}

\begin{pf} We get {\it i)} by definition and  {\it ii)} is just the
Corollary above. To prove {\it iii)},
let us write (\ref{funrel}) from Lemma 1 in the form
$$\tau_{N^i}D_{N^k}\mu=D_{N^{k+i}}\mu- \delta_{N^{i}}D_{N^k}\mu\ .
$$
Hence,
$$D_{N^{k+i}}\mu=\delta_{N^{i}}D_{N^k}\mu+\tau_{N^i}D_{N^k}\mu=D_{N^{i}}
D_{N^{k}}\mu\ .$$
Finally,  {\it iv)} is straightforward. By the result of Lemma 1,
\bea
N^i(D_{N^{k+i}}\mu(X,Y))&=&(A_N^iD_{N^{k+i}}\mu)(X,Y)=(A_N^i\delta_{N^{i
}}+T_{N^i})
D_{N^k}\mu(X,Y)
\cr &=&
B_N^iC_N^iD_{N^k}\mu(X,Y)=D_{N^k}\mu(N^i(X),N^i(Y))\ .
\eea
\end{pf}
\begin{cor} For any strong Saletan tensor $N$ for $\zm$,

i) $N^k({\cal A}) $ is a subalgebra with respect to the product
$D_{N^i}\mu$;

ii)  $\ker  N^k=\{X\in{\cal  A}\mid  N^k(X)=0  \}$  is   an   ideal   of
$D_{N^i}\mu$, for all $i>k$.
\end{cor}

\section{Behaviour  of  properties   of   algebraic   structures   under
contraction}

Assume that our product $\zm$ is a specific one, satisfying some general
axioms $\{ (a_\zm^i)\}$ of the form
\be\label{axioms}
(a^i_\zm)\qquad\forall x_1,\dots,x_{n_i}\in\A\quad
[w_\zm^i(x_1,\dots,x_{n_i})=0],
\ee
where $w_\zm^i$ are $\zm$-polynomial functions, and using only universal
quantifiers, like

$$\leqno{(a^1_\zm)}
\qquad\forall x,y,z,\in\A \qquad
[\mu(x,\mu(y,z))+\mu(y,\mu(z,x))+\mu(z,\mu(x,y))=0],
$$
or
$$\leqno{(a^2_\zm)}
\qquad\forall x,y\in\A \qquad [\mu(x,y)+\mu(y,x)=0]\ ,
$$
or
$$\leqno{(a^3_\zm)}
\qquad\forall x,y,z\in\A\qquad [\mu(x,\mu(y,z))-\mu(\mu(x,y),z)=0]\ ,
$$
but not using existential quantifiers like
$$\leqno{(a^4_\zm)}
\qquad\exists 1\in\A\ \forall y\in\A\quad [\mu(1,y)=y=\mu(y,1)]\,.
$$
An algebra satisfying $(a^1_\zm)$ and $(a^2_\zm)$ is a Lie algebra,
an algebra satisfying $(a^3_\zm)$ is associative, and $(a^4_\zm)$
says that $\A$ is unital.

\begin{theorem} If the product $\zm$ satisfies axioms of the form
(\ref{axioms}), then the  contracted  product  $\zm_N$  satisfies  these
axioms.
\end{theorem}
\begin{pf}     The     products      $\zm^\zl_N=U(\zl)^{-1}\circ\zm\circ
U(\zl)^{\otimes 2}$ are isomorphic to $\zm$, so that  they  satisfy  the
same  axioms, and  equations  $w_{\zm^\zl_N}^i(x_1,\dots,x_{n_i})=0$   are
going to $w_{\zm_N}^i(x_1,\dots,x_{n_i})=0$ by passing to the limit as
$\zl\to 0$.
\end{pf}

\medskip\noindent
{\bf Remark.} The above theorem implies that  a  contraction  of  a  Lie
algebra is a Lie algebra and a contraction of an associative algebra  is
an associative algebra. However, it is crucial that the axioms use the
universal quantifiers only. For example, the existence of unity
($a^4_\zm$) is, in general, not preserved by contractions as shows the
case $N=0$.  This  is because the unit for the product $\zm_N^\zl$ is
$U(\zl)^{-1}(1)$ which may have no limit as $\zl\to 0$.

\medskip\noindent
{\bf Definition.} We say that products $\zm,\zm'$ satisfying axioms
(\ref{axioms})  are  {\it  compatible},  if   any   linear   combination
$\zm+\za\zm'$ satisfies these  axioms.  For  instance,  two  associative
products are compatible if and only if its sum is associative  as  well,
etc.

\begin{theorem} If $N$  is  a  Nijenhuis  tensor  for  $\zm$,  then  the
products $\zm$ and $\zm_N=\zd_N\zm$ are compatible.
\end{theorem}
\begin{pf}
According to Theorem 3, $I+\za\, N$ is a Nijenhuis tensor for  $\zm$  for
any $\za\in\R$. Using now Theorem 10 we see that the product
$$\zm_{(I+\za N)}=\zd_{(I+\za N)}\zm=\zm+\za\,\zm_N
$$
satisfies the axioms of $\zm$.
\end{pf}

\medskip\noindent
{\bf Remark.} If $N$ is only a Saletan tensor, the  products  $\zm$  and
$\zm_N$ are, in general, not  compatible.  For  example,  the associative
products $X*Y$ and $X*_NY=N^{-1}(N(X)*N(Y))$, for invertible $N$, are  in
general not compatible, i.e., $X*Y+X*_NY$ is, in general, not associative.

\section{Contractions of Lie algebras}

Let us consider now the  very important particular case of a
finite-dimensional Lie algebra $(E ,[\cdot,\cdot])$. This corresponds to
the vector bundle $E$ over a single point with $\A=E$ and
$\zm=[\cdot,\cdot]$.
The family $U(\ze)$ of endomorphisms of the underlying vector space $V$
considered by \.In\"on\"u and Wigner \cite{IW} is
$U(\ze)=P+\ze (I-P)$, where $P$ is a projection,  and it was
later studied by Saletan \cite{Sa} in the more general case
$$
U(\epsilon)=\epsilon{\,I}+(1-\epsilon){N},
$$
for which $U(0)=N$ and $U(1)=I$. By reparametrizing it  with a
new parameter $\lambda=\frac{\epsilon}{1-\epsilon} $,
it is, as shows Theorem 4, equivalent to the contraction
\be
[X,Y]_N=\lim_{\zl\to 0}U(\lambda)^{-1}
[U(\lambda)X,U(\lambda)Y],
\ee
with $U(\zl)=\zl\, I+N$.
In this particular case, the Riesz decomposition $E=E^1\oplus E^2$ with
respect to $N$ is regular and, according to our general Theorem 2 and
Theorem 3, the necessary and
sufficient condition for the existence of such a limit is that
\begin{equation}
T_N\zm(X,Y)=[NX,NY]-N[NX,Y]-N[X,NY]+{N^2}[X,Y]\in E^1\, . \label{eq:CNS}
\end{equation}
Moreover, we obtain the following expression for the new bracket:
$$
[X,Y]_N=\left.
N\right|_{E_1}^{-1}[NX,NY]_1-N[X,Y]_2+[NX,Y]_2+[X,NY]_2,
$$
where the subscripts refer to the projections onto $E^1$ or $E^2$.
Consequently, Theorem 2 implies
\begin{eqnarray}
N[X,Y]_N=[NX,NY]\ .
\end{eqnarray}
Therefore, a necessary  condition for the existence of a contraction
leading from a Lie algebra to become another one is the existence of a Lie
algebra homomorphism of the second into the first one. However, as
Levy--Nahas pointed out this is not a sufficient condition \cite{LN}.

The necessary and sufficient condition as expressed by Gilmore in
\cite{G}:

\smallskip\noindent
{\sl the contraction exists if and only if}
\begin{equation}
{N^{p+s}[X,Y]_2-N^{p}[X,N^{s}Y]_2}={N^{s}[N^{p}X,Y]_2-[N^{p}X,N^{s}Y]_2}
\quad {\rm for\ all\ }p,s>0\,,
\label{Gilmore}
\end{equation}
can be easily obtained using techniques developed in Section 4. Indeed, in
the notation of Section 4, (\ref{Gilmore}) reads
\be\label{gi}
(A_N^{p+s}-A_N^pC_N^s)\zm_2=(A_N^sB_N^p-B_N^pC_N^s)\zm_2,
\ee
where $\zm_2$ is the projection of the bracket onto $E_2$. Since all
operators commute among themselves and with the projection, we can write
(\ref{gi}) in the form
\beas
(A_N^s-C_N^s)(A_N^p-B_N^p)\zm_2&=&w(A_N,B_N,C_N)
((A_N-C_N)(A_N-B_N)\zm)_2 \\
&=&w(A_N,B_N,C_N)(T_N\zm)_2=0,
\eeas
which is true for all $p,s>0$ if and only if $(T_N\zm)_2=0$, since the
polynomial $w$ equals $1$ for $p=s=1$.

\medskip\noindent
{\bf Example 5.} Using Theorem 5 we get the \.In\"on\"u-Wigner contraction
for Lie algebras. Consider just a splitting $E=E^1\oplus E^2$ of the Lie
algebra $E$ into a subalgebra $E^1$ and a complementary subspace $E^2$.
According to Theorem 5, the projection $N$ of $E$ onto $E^1$ along $E^2$
is a Saletan tensor with the splitting being also the Riesz decomposition.
The resulting bracket is
\be\label{81}
[X,Y]_N=[X_1,Y_1]+[X_1,Y_2]_2+[X_2,Y_1]_2.
\ee
To have a particular example, take $E=\su$ with the basis $X_1,X_2,X_3$
satisfying the commutation rules
\be\label{82}
[X_1,X_2]=X_3,\quad [X_2,X_3]=X_1,\quad [X_3,X_1]=X_2.
\ee
As for $E^1$, take the 1-dimensional subalgebra spanned by $X_1$, and let
$E^2$ be spanned by $X_2,X_3$. According to (\ref{81}), the commutation
rules for the contracted algebra read
\be\label{83}
[X_1,X_2]=X_3,\quad [X_2,X_3]=0,\quad [X_3,X_1]=X_2.
\ee
One recognizes easily the Lie algebra ${\goth e}(2)$ of Euclidean motions in
in a two-dimensional space.

\medskip
As there were some {contractions} that could not be explained neither
in the framework of \.In\"on\"u--Wigner \cite{IW} nor in that of Saletan \cite{Sa},
Levy--Nahas proposed a more singular contraction procedure  by assuming
families $U(\zl)={\zl^p}\,{U_{s}(\zl)}$, where $p\in{\Bbb N}$ and
$U_{s}(\zl)=N+\zl\, I$.
Following a quite similar path to that of Saletan contractions,
one obtains as a necessary and sufficient condition for the existence of
the limit in the $p=1$ case,
\begin{equation}
N(T_N(X,Y)_2)=0  \ ,  \label{s11}
\end{equation}
where
$$
T_N(X,Y)_2=[NX,NY]_2-N[X,NY]_2-N[NX,Y]_2+{N^2}[X,Y]_2
$$
is the projection of the Nijenhuis torsion onto $E^2$ (cf. with
the condition $T_N(X,Y)_2=0$ in the standard case).
The new bracket is then $T_N(X,Y)_2$
For general $p$, the condition for the existence of the contraction is
$N^p(T_N(X,Y)_2)=0$, and the resulted bracket is $(-N)^{p-1}T_N(X,Y)_2$.

For the sake of completeness we will finally mention that  other generalized
\.In\"on\"u-Wigner contractions were proposed in  \cite{DM} and
\cite{WW1}.

\section{Contractions of Lie algebroids}

Lie algebroids, which are very common structures in geometry, should be
very nice objects for contractions in our sense, since they are, by
definition, certain algebra structures on sections of vector bundles. They
were introduced by Pradines \cite{Pr} as infinitesimal objects for
differentiable groupoids, but one can find similar notions proposed by
several authors in increasing number of papers (which proves their
importance and naturalness). For basic properties and the literature on the
subject we refer to the survey article by Mackenzie \cite{Mac}.
\begin{definition} A {\sl Lie algebroid} on a smooth manifold $M$ is a
vector bundle $\zt:E\ra M$, together with a bracket
$\zm=[\cdot,\cdot]:\A\ti\A\ra\A$ on the $C^\infty(M)$-module $\A=\zG(E)$
of smooth sections of $\zt$, and a vector bundle morphism $a_\zm:E\ra TM$,
over the identity on $M$, from $E$ to the tangent bundle $TM$, called the
{\sl anchor} of the Lie algebroid, such that

(i) the bracket $\zm$ is a Lie algbera bracket on $\A$ over $\R$;

(ii) for all $X,Y\in\A$ and all smooth functions $f$ on $M$ we have
\be\label{la1}
\zm(X,f\,Y)=f\,\zm(X,Y)+a_\zm(X)(f)\,Y;
\ee

(iii) For all $X,Y\in\A$,
\be\label{la2}
a_\zm(\zm(X,Y))=[a_\zm(X),a_\zm(Y)],
\ee
where the square bracket is the Lie bracket of vector fields. In other words,
$a_\mu$ is a Lie algebra homorphism.
\end{definition}
\medskip\noindent
{\bf Example 6.} Every finite-dimensional Lie algebra $E$ is a Lie
algebroid
as a bundle over a single point with the trivial anchor. More generally,
any bundle of Lie algebras is a Lie algebroid with the trivial anchor.

\medskip\noindent
{\bf Example 7.} There is a canonical Lie algebroid structure on every
tangent bundle $TM$ with the bracket being the standard bracket of vector
fields and the anchor being just the identity map on $TM$.

\medskip\noindent
{\bf Example 8.} There is a natural Lie algebroid associated with a
realization of a Lie algebra in terms of vector fields. Suppose $V$ is a
Lie algebra with the bracket $[\cdot,\cdot]$ with a realization
$\hat{} :V\ra{\goth X}(M)$ in terms of vector fields on a manifold $M$.
We can view $V$ as a subspace of sections
on the trivial bundle $E=M\ti V$ over $M$, regarding $X\in V$ as constant
sections of $E$. There is uniquely defined Lie algebroid structure on
$\A=\zG(E)=C^\infty(M,V)$
such that the Lie algebroid bracket $\zm$ and the anchor $a_\zm$ satisfy:

(i) $\zm(X,Y)=[X,Y]$ for all $X,Y\in V$;

(ii) $a_\zm(X)=\hat X$ for every $X\in V$.

\noindent
In other words, identifying $\A$ with $C^\infty(M)\otimes V$, the Lie
algebroid bracket reads
\be
\zm(f\otimes X,g\otimes Y)=fg\otimes[X,Y]+f\hat X(g)\otimes Y
-g\hat Y(f)\otimes X.
\ee

\medskip
\noindent{\bf Example 9.} There is a canonical Lie algebroid structure on the
cotangent bundle $T^*M$ associated with a Poisson tensor $P$ on $M$.
This is the unique Lie algebroid bracket $[\cdot,\cdot]^P$ of differential
1-forms for which $[\D f,\D g]^P=\D\{ f,g\}^P$, where
$\{\cdot,\cdot\}^P$ is the Poisson bracket of functions for $P$, and the
anchor map is just $P$ viewed as a bundle morphism $P:T^*M\ra TM$.
Explicitly,
\be
[\za,\zb]^P=\Li_{P(\za)}\zb-\Li_{P(\zb)}\za-\D\<P,\za\we\zb>.
\ee
This Lie bracket was defined first by Fuchssteiner \cite{Fu}.
We shall comment more on this structure in the next section.

\medskip
It is interesting that any contraction of a Lie algebroid bracket
gives again a Lie algebroid bracket. We start with the following lemma.
\begin{lemma} If $\zm$ is a Lie algebroid bracket on $\A=\zG(E)$ and
$a_\zm:E\ra TM$ is the corresponding anchor, then for any $(1,1)$--tensor
$N$ on $E$ we have

(i) $\zd_N\zm(X,f\,Y)=f\,\zd_N(X,Y)+a_\zm(N(X))(f)\,Y$;

(ii) $T_N\zm(X,f\,Y)=f\,T_N\zm(X,Y)$,

\smallskip\noindent
for any $X,Y\in\A$, $f\in C^\infty(M)$.
\end{lemma}
\begin{pf} {\it (i)} By definition and properties of Lie algebroid brackets,
\beas
\zd_N\zm(X,f\,Y)&=&\zm(N(X),f\,Y)+\zm(X,N(f\,Y))-N\zm(X,f\,Y)\\
&=&f\,\zm(N(X),Y)+a_\zm(N(X))(f)\,Y+f\,\zm(X,N(Y))\\
&&+a_\zm(X)(f)\,N(Y)-f\,N\zm(X,Y)+a_\zm(X)(f)\,N(Y)\\
&=&f\,(\zm(N(X),Y)+\zm(X,N(Y))-N\zm(X,Y))+a_\zm(N(X))(f)Y\\
&=&f\,\zd_N\zm(X,Y)+a_\zm(N(X))(f)\,Y.
\eeas
Here we have used the fact that the multiplication by a function commutes
with $N$ (i.e., $N$ is a tensor).

{\it (ii)}\ We have
\beas
T_N\zm(X,f\,Y)&=&\zm(N(X),N(f\,Y))-N\zd_N\zm(X,f\,Y)\\
&=&f\,\zm(N(X),N(Y))+a_\zm(N(X))(f)\,N(Y)\\
&&-f\,N\zd_N\zm(X,Y)-a_\zm(N(X))(f)\,N(Y)\\
&=&f\,(\zm(N(X),N(Y))-\,N\zd_N\zm(X,Y))=fT_N\zm(X,Y),
\eeas
where we have used (i).
\end{pf}
\begin{theorem} If $N$ is a Saletan tensor for a Lie algebroid bracket
$\zm$ on $\A=\zG(E)$, with an anchor map $a_\zm:E\ra TM$, then the
contracted bracket $\zm_N$ is again a Lie algebroid bracket on $\A$ with
the anchor $a_{\zm_N}=a_\zm\circ N$.
\end{theorem}
\begin{pf} We already know that the contracted bracket
$\zm_N=\zd_N\zm+\zt_N\zm$ is a Lie bracket.
Since $N$ is a Saletan tensor, $\zt_N\zm=N^{-1}T_N\zm$ is well-defined and
clearly satisfies also (ii) of the above Lemma. Thus, using also (i),
\beas
\zm_N(X,fY)&=&\zd_N\zm(X,fY)+\zt_N\zm(X,fY)\\
&=&f\zd_N\zm(X,Y)+a_\zm(N(X))(f)Y+f\zt_N\zm(X,Y)\\
&=&f(\zd_N\zm(X,Y)+\zt_N\zm(X,Y))+a_\zm(N(X))(f)Y\\
&=&f\zm_N(X,Y)+a_\zm(N(X))(f)Y,
\eeas
so that $a_{\zm_N}=a_\zm\circ N$ can serve for the anchor of $\zm_N$.
It suffices to check the condition (\ref{la2}):
\beas
[a_{\zm_N}(X),a_{\zm_N}(Y)]&=&[a_\zm(N(X)),a_\zm(N(Y))]\\
&=&a_\zm(\zm(N(X),N(Y))=a_\zm(N\zm_N(X,Y)).\\
\eeas
We have used the identity $\zm(N(X),N(Y))=N\zm_N(X,Y)$ which holds for
Saletan tensors.
\end{pf}

Note that this type of contractions of Lie algebroids has been already
studied by Kosmann-Schwarzbach and Magri in \cite{KSM} in the case of
Nijenhuis tensors. All results of this section can also be applied to
general algebroids as defined in \cite{GU}.

\medskip\noindent
{\bf Example 10.} The contracted bracket of vector fields defined in
Example 1 defines a new Lie algebroid structure on $TM$ with the anchor
map being the projection onto the subbundle $E^1$ of $TM$.

\medskip\noindent
{\bf Example 11.} Any Saletan contraction of a Lie algebra $V$ leads to a
contraction of the Lie algebroid associated with an action of $V$ on $M$,
which was described in Example 8. More precisely, if $N_0$ is a Saletan
tensor for $V$, then
\be
N:M\ti V\ra M\ti V, \quad N(f\otimes X)=f\otimes N_0(X),
\ee
is a Saletan tensor for the canonical Lie algebroid bracket on
$\A=C^\infty(M)\otimes V$. Indeed, if $V=V^1\oplus V^2$ is the Riesz
decomposition for $N_0$, then $E=E^1\oplus E^2$, with
$E^i=C^\infty(M)\otimes V^i$, is the Riesz decomposition for $N$.
Moreover,
the Nijenhuis torsion $T_N\zm$ takes values in $\A^1$. Indeed,
by Lemma 2(ii), the Nijenhuis torsion $T_N\zm$ is tensorial, so it
suffices to check that on $V$ it takes values in $\A^1$. But on $V$ the
Nijenhuis torsion of $N$ with respect to $\zm$ is the same as the
Nijenhuis torsion of $N_0$ with respect to the bracket on $V$, so it takes
values in $V^1\subset\A^1$. Finally, $E^1$ is of constant rank, so $N$ is
regular and hence a Saletan tensor due to Theorem 3. The anchor map for
$\zm_N$ is $a_\zm\circ N$, so $a_{\zm_N}(f\otimes X)=f\widehat{N_0(X)}$
and the contracted anchor takes values in the module of vector fields
generated by the action of the subalgebra $N_0(V)$ on $M$. In fact,
what we get is the Lie algebroid structure on $M\ti V$ associated with the
contracted Lie algebra structure on $V$ and the anchor map $a_\zm\circ
N$.

As a particular example let us take the Lie algebroid on $S^2\ti\su$
associated with the action of the Lie algebra $\su$ on the 2-dimensional
sphere $S^2=\{(x,y,z)\in\R^3:x^2+y^2+z^2=1\}$ given by (in the notation of
Example 5)
\be
\hat X_1=y\pa_z-z\pa_y,\quad \hat X_2=x\pa_z-z\pa_x,\quad \hat
X_3=x\pa_y-y\pa_x.
\ee
From the contraction of $\su$ into ${\goth e}(2)$, as described in Example
5, we construct a contraction of this Lie algebroid. For the Lie bracket
we get
\beas
\zm_N\left(\sum f_i\otimes X_i,\sum g_j\otimes X_j\right)
&=&(f_1\hat X_1(g_1)-g_1\hat X_1(f_1))\otimes X_1\\
&&+(f_3g_1-f_1g_3+f_1\hat X_1(g_2)-g_1\hat X_1(f_2))\otimes X_2\\
&&+(f_1g_2-f_2g_1+f_1\hat X_1(g_3)-g_1\hat X_1(f_3))\otimes X_3
\eeas
and the anchor reads
\be
a_{\zm_N}\left(\sum_{i=1}^3f_i\otimes X_i\right)=f_1\hat X_1.
\ee
This is the Lie algebroid structure on $S^2\ti{{\goth e}}(2)$ associated with
the representation $\widehat{(X_i)_N}=\zd^i_1\hat X_1$ of ${{\goth e}}(2)$
 in terms of vector fields on $S^2$.

\section{Poisson contractions}

Poisson brackets, being defined on functions, are brackets of sections of
1-dimensional bundles and seem, at the first sight, not to go under our
contraction procedures. We shall show that this is not true and that our
contraction method shows precisely what contraction of a Poisson tensor
should be. The crucial point is that we should think about a Poisson
tensor $P$ on $M$ as defining certain Lie algebroid structure on $T^*M$
rather than defining just the Poisson bracket $\{\cdot,\cdot\}^P$ on
functions.
Recall from Example 9 that the Lie algebroid bracket on differential
forms, associated with $P$, reads
\be\label{la}
[\za,\zb]^P=\Li_{P(\za)}\zb-\Li_{P(\zb)}\za-\D\<P,\za\we\zb>.
\ee
The anchor map of this Lie algebroid is just $P$, viewed as a bundle
morphism $P:T^*M\ra TM$, so that $P$ can be easily decoded from the Lie
algebroid structure. Of course, not all Lie algebroid structures on
$T^*M$, even
having a Poisson tensor for the anchor map, are of this kind. We just
mention that an elegant characterization of Lie algebroid brackets
associated with Poisson structures is that the exterior derivative
acts as a graded derivative on the corresponding Schouten-like bracket of
differential forms, or that this Lie algebroid structure constitutes a
{\it Lie bialgebroid} (in the sense of Mackenzie and Xu \cite{MX})
together with the canonical Lie algebroid structure on the dual, i.e.,
the tangent bundle $TM$.

One can think differently. Suppose we have a skew-symmetric 2-vector field
$P$, viewed as a bundle morphism $P:T^*M\ra TM$, and we write formally
the bracket (\ref{la}). When do we obtain a Lie algebroid bracket? The
answer is very simple (cf.\cite{KSM}):
\begin{theorem} If $P$ is a skew 2-vector field, then formula
(\ref{la}) gives a Lie algebroid bracket if and only if $P$ is a Poisson
tensor.
\end{theorem}

Let us fix now a Poisson structure $P$ on $M$ and the Lie algebroid
bracket $\zm^P=[\cdot,\cdot]^P$. Given Saletan tensor $N$ for $\zm^P$, we
get the contracted bracket $\zm^P_N$. It is natural, in the case when
$\zm^P_N$ is again a Lie algebroid bracket associated with a Poisson
tensor (we shall speak about {\it Poisson contraction}), to call this
tensor a {\it contracted Poisson structure} by means of $N$.
By $^tN$ we shall denote the dual bundle morphism of $N:T^*M\ra T^*M$. In
particular, $^tN:TM\ra TM$. We almost follow the notation of \cite{MM,KSM}
but with exchanged roles for $N$ and $^tN$ which, as it will be seen
later, seems to be more appropriate in our case.

\begin{theorem}
A necessary and sufficient condition for the contraction of the Lie
algebroid $\zm^P$ associated with a Saletan tensor N to be a Poisson
contraction is that

(i) $PN={^tN}P$

and

(ii) $\zm^P_N=\zm^{PN}$.
\end{theorem}
\begin{pf} First, assume that the contraction according to $N$ is a
Poisson contraction. Hence, $\zm^P_N=\zm^{P_1}$ for a Poisson tensor
$P_1$. But the anchor of $\zm^{P_1}$ is $P_1$ and the anchor of $\zm^P_N$
is $PN$ (Theorem 12). We get then $P_1=PN$ and (ii). Since $PN$ must be
skew-symmetric, $^t(PN)=-PN$. But $^t(PN)={^tN}{^tP}=-{^tN}P$ and we get
(i).

Suppose now (i) and (ii). Since (i) means that $PN$ is skew-symmetric and
$\zm^P_N$ is a Lie algebroid bracket, in view of Theorem 13, the tensor
$PN$ is a Poisson tensor.
\end{pf}

\medskip\noindent
{\bf Remark.} Bihamiltonian systems, as noticed by Magri \cite{Mag}, play
an important role in the discussion of complete integrability in the sense
of Liouville. A geometrical approach to this questions, proposed in
\cite{MM} (see also \cite{KSM}), uses the notion of a {\it
Poisson-Nijenhuis structure}, i.e., a pair $(P,N)$, where $P$ is a Poisson
tensor on $M$ and $N$ is a Nijenhuis tensor on the tangent bundle $TM$,
which satisfy certain compatibility conditions. For contractions of
$\zm^P$, we use $N$ being a morphism of $T^*M$ rather than of $TM$, but
of
course, by duality,  $^tN:TM\ra TM$.

In the case when $N$ is a Nijenhuis tensor for $\zm^P$, our conditions
(i) and (ii) are the same as the
compatibility conditions for Poisson-Nijenhuis structure of \cite{MM,KSM}
with $N$ replaced by $^tN$. In this case  $\zm^{PN}-\zm^P_N$
is exactly what in \cite{KSM} is denoted by by $C(P,{^tN})$.
Note that Poisson-Nijenhuis structures can be described in terms of Lie
bialgebroids \cite{KS} (cf. also \cite{GU1} for a more general setting).

We do not assume that $^tN$ (in
our notation) is a Nijenhuis tensor for the canonical Lie algebroid $TM$,
but that $N$ is a Saletan tensor for $\zm^P$ on $T^*M$.
It is natural to call a pair $(P,N)$, where $P$ is a Poisson structure on
$M$ and $N$ is a Saletan tensor for $\zm^P$ satisfying (i) and (ii) of the
above theorem, a {\it Poisson-Saletan structure}.
If $(P,N)$ is a Poisson-Saletan structure, then $(P,{^tN})$ need not be a
Poisson-Nijenhuis structure in the sense of \cite{KSM},
even if we impose that $N$ is a Nijenhuis tensor for $\zm^P$,
as shows the following example. However, this weaker assumption is
sufficient to perform a Poisson contraction and to obtain the contracted
Poisson structure $PN$. In the case when $N$ is a Nijenhuis tensor for
$\zm^P$, according to Theorem 11, Lie algebroid brackets $\zm^P$ and
$\zm^P_N=\zm^{PN}$ are compatible, so also the Poisson tensors $P$ and
$PN$ are compatible. We can also get a whole hierarchy of compatible
Poisson tensors using the results of Section 4.

\medskip\noindent
{\bf Example 12.} Let $M=M_1\ti M_2$, where $M_i$, $i=1,2$, is a manifold.
On the product manifold consider the product Poisson structure $P=P_1\ti
\{ 0\}$, where $P_1$ is a Poisson structure on $M_1$. Let
$N_2:T^*M_2\ra T^*M_2$ be any $(1-1)$--tensor on $M_2$. It induces a
tensor $N:T^*M\ra T^*M$ which on
\be
T^*_{(m_1,m_2)}M=T^*_{m_1}M_1\oplus T^*_{m_2}M_2
\ee
acts by identity on $T^*_{m_1}M_1$ and by $(N_2)_{m_2}$ on
$T^*_{m_2}M_2$. The $C^\infty(M)$-module $\zW^1(M)$ of 1-forms on $M$ is
generated by $\zW^1(M_1)$ and $\zW^1(M_2)$ and, as can be easily seen from
(\ref{la}),
$\zm^P(\za,\zb)=\zm^{P_1}(\za,\zb)$ for $\za,\zb\in\zW^1(M_1)$, and
$\zm^P(\za,\zb)=0$ when $\za\in\zW^1(M_2)$. Since
$\zW^1(M_1)$ and $\zW^1(M_2)$ are invariant subspaces for $N$, and
since $N$ acts by identity on $\zW^1(M_1)$, it follows that $N$ is a
Nijenhuis tensor for $\zm^P$ and that $\zm^P_N=\zm^P=\zm^{PN}$. Thus
$(P,N)$ is a
Poisson-Saletan structure. On the other hand, $^tN=id\ti{^tN_2}$ need not
to be a Nijenhuis tensor for $TM$, since $N_2$ is arbitrary.

\medskip\noindent
However, we have the following weaker result.

\begin{theorem} If $N:T^*M\ra T^*M$ is a Nijenhuis tensor for $\zm^P$,
then the Nijenhuis torsion of $^tN:TM\ra TM$ vanishes on the vector
fields from the image of $P:T^*M\ra TM$. In particular, if $P$ is
invertible, i.e., comes from a symplectic structure, then $(P,{^tN})$ is
a Poisson-Nijenhuis structure is the sense of \cite{KSM}.
\end{theorem}
\begin{pf} Writing down $T_N\zm^P=0$, we get
\be
\zm^P(N(\za),N(\zb))=N(\zm^P(N(\za),\zb)+\zm^P(\za,N(\zb))-N\zm^P(\za,\zb)).
\ee
Applying the anchor $P$ to both sides, we get, according to (\ref{la2}),
\be
[PN(\za),PN(\zb)]=PN(\zm^P(N(\za),\zb)+\zm^P(\za,N(\zb))-N\zm^P(\za,\zb)).
\ee
Using now $PN={^tN}P$ and the fact that anchor is a homomorphism of
the brackets once more, we get
\be
[{^tN}P(\za),{^tN}P(\zb)]={^tN}([{^tN}P(\za),P(\zb)]+[P(\za),{^tN}P(\zb)]
-{^tN}[P(\za),P(\zb)]).
\ee
The last means exactly that
\be\label{week}
T_{^tN}(P(\za),P(\zb))=0,
\ee
where $T_{^tN}$ is the Nijenhuis torsion of $^tN$ with respect to the
bracket of vector fields.
\end{pf}
{\bf Remark.} The property (\ref{week}), together with the compatibility
condition, defines a {\it week Poisson-Nijenhuis structure} in the
terminology of \cite{MMP}. That this week condition is sufficient to get
recursion operators was first observed in \cite{MN}. Note also that a
similar procedure can be applied to Jacobi structures. Jacobi structures
give rise to Lie algebroids as was observed in \cite{KSB}. Similarly as
above, the contraction procedures for these Lie algebroids give rise to a
proper concept of a {\it Jacobi-Nijenhuis structure}. We refer to
\cite{MMP} for details.

\section{Contractions of $n$-ary products and
coproducts}

Let, as before,  $N$  be  a  (1,1)-tensor  over  a  vector  bundle  $E$,
$E=E^1\oplus E^2$ be the Riesz decomposition of $E$ relative to $N$, and
$\A,\A^1$ be the spaces of smooth sections of $E$ and $E^1$, respectively.
In complete analogy with binary  products, we  can  consider  $n$-ary
products (resp. coproducts),  i.e.,  linear  mappings
$\zm:\A^{\otimes n}\ra\A$ (resp. linear mappings $\zm:\A\ra\A^{\otimes
n}$) and contractions of them with  respect  to  families  $U(\zl)=\zl\,
I+ N$. For an $n$-ary product (resp. coproduct) we denote
$$\zd_N\zm=\zm\circ N^{n-1}_n-N\circ\zm\circ N^{n-2}_n+\dots +(-1)^{n-1}
N^{n-1}\circ\zm\circ N^0_{n},
$$
and, respectively,
$$\zd_N\zm=N^{n-1}_n\circ\zm-N^{n-2}_n\circ\zm\circ N+\dots
+(-1)^{n-1} N^0_n\circ\zm\circ N^{n-1},
$$
where $N^k_n$ are defined by
$$(\zl\, I+N)^{\otimes n}=\sum_{k=0}^n\zl^{n-k}N^k_n.
$$
An obvious adaptation of the proof of Theorem 2 gives the following.
\begin{theorem}
Let  $\mu:{\cal  A}^{\otimes n}\ra{\cal  A}$   be   a
point-wise continuous $n$-ary product in $\A$.
Then, for $U(\zl)=\zl\, I+N$, the limit
$$\zm_N=\lim_{\lambda\to 0}(U(\lambda)^{-1}\circ\zm\circ U(\zl)^{\otimes
n})
$$
exists and defines a new (contracted) $n$-ary product
$D_N\zm$ on $\cal A$ if and only if the Nijenhuis torsion
$$T_N\mu=\zm\circ N^{\otimes n}-N\circ\zd_N\zm
$$
takes values in $N({\cal A}^1)$. If this is the case, then
\be
\zm_N=\zd_N\zm+\tau_N\zm,
\ee
where $N(\tau_N\mu)=T_N\mu$.
Moreover, $N$ constitutes a homomorphism of $({\cal A},\mu_N)$ into
$({\cal A},\mu)$:
$$N\circ\zm_N=\zm\circ N^{\otimes n}\ .
$$
\end{theorem}

A similar theorem for coproducts can be obtained by duality. Since it is
much harder to put conditions for existence of contraction at  arguments
of the Nijenhuis torsion, for simplicity we give an explicit version for
the regular case only.

\begin{theorem}
Let  $\mu:{\cal  A}\ra{\cal  A}^{\otimes n}$   be   a
point-wise continuous $n$-ary coproduct in $\A$. If $N$ is regular, i.e.,
$E^1$ is of constant dimension, then, for $U(\zl)=\zl\, I+N$, the limit
$$\zm_N=\lim_{\lambda\to  0}(U(\lambda)^{\otimes n}\circ\zm\circ
U(\zl)^{-1})
$$
exists and defines a new (contracted) $n$-ary coproduct
$D_N\zm$ on $\cal A$ if and only if the Nijenhuis torsion
$$T_N\mu=N^{\otimes n}\circ\zm-\zd_N\zm\circ N
$$
vanishes on ${\cal A}^2$-the space of sections of $E^2$. If this is the
case, then
\be
\zm_N=\zd_N\zm+\tau_N\zm,
\ee
where $(\tau_N\mu)=T_N\mu\circ N^{-1}$ on $\A^1$  and  $\tau_N\zm=0$  on
$\A^2$.
Moreover, $N$ constitutes a homomorphism of $({\cal A},\mu_N)$ into
$({\cal A},\mu)$:
$$\zm_N\circ N=N^{\otimes n}\circ\zm \ .
$$
\end{theorem}

One  can  consider  more  general  algebraic  structures of   the   form
$\zm:\A^{\otimes k}\ra\A^{\otimes n}$, but this leads to more conditions
of contractibility and we will not study  these  cases  in  the  present
paper. Note only that contractions of coproducts, as a part of
contractions of Lie bialgebras, appeared already in \cite{BGHO}.

\section{Conclusions}

Motivated by physical examples from Quantum Mechanics, we have studied
contractions of general binary (or $n$-ary) products with respect to
one-parameter families of transformations  of  the  form  $U(\zl)=\zl\,
A+N$, generalizing pioneering work by \.In\"on\"u, Wigner, and Saletan.
Our generalization can be applied to many infinite-dimensional cases,
especially Lie algebroids and Poisson brackets, however, it does not deal
with a generic dependence of the contraction parameter. The problem of
describing  contractions  with  respect  to   general   $U(\zl)$, or even
differentiable with respect to $\zl$, is much  more  complicated.

The contraction procedure can be viewed as an inverse of a  deformation
procedure. Deformations of associative and Lie algebras, at least on the
infinitesimal level,  are  related  to  some  cohomology.  It  would  be
interesting to relate formally deformations to contractions  and  connect
the cohomology also to contractions.

We postpone these problems to a separate paper.


\end{document}